\documentclass[ECP]{null} 
\SHORTTITLE{Enhanced Stein Couplings}

\TITLE{Zero Bias Enhanced Stein Couplings}

\AUTHORS{
  Larry~Goldstein\footnote{Department of Mathematics, University of Southern California, Los Angeles, CA USA
    \EMAIL{larry@math.usc.edu}}}


\AMSSUBJ{60F05; 62P10; 60C05} 

\SUBMITTED{June 10, 2022} 
\ACCEPTED{00.00.00}

\VOLUME{0}
\YEAR{2020}
\PAPERNUM{0}
\DOI{10.1214/YY-TN}

\ABSTRACT{The Stein couplings of Chen and Roellin \cite{chen2010stein} vastly expanded the range of applications for which coupling constructions in Stein's method for normal approximation could be applied, and subsumed both Stein's classical exchangeable pair, as well as the size bias coupling. A further simple generalization includes zero bias couplings, and also allows for situations where the coupling is not exact. The zero bias versions result in bounds for which often tedious computations of a variance of a conditional expectation is not required. An example to the Lightbulb process shows that even though the method may be simple to apply, it may yield improvements over previous results that had achieved bounds with optimal rates and small, explicit constants.}

\usepackage{bm,xcolor}
\newcommand{\qmq}[1]{\quad \mbox{#1} \quad}
\newcommand{\qm}[1]{\quad \mbox{#1}}
\newcommand{\bbox}{\hfill $\Box$}
\newcommand{\E}{\mathbb{E}}
\newcommand{\ycolor}[1]{\textcolor{blue}{#1}}
\newcommand{\bea}{\begin{eqnarray}}
	\newcommand{\ena}{\end{eqnarray}}
\newcommand{\beas}{\begin{eqnarray*}}
	\newcommand{\enas}{\end{eqnarray*}}

\providecommand{\bysame}{\leavevmode\hbox to3em{\hrulefill}\thinspace}
\providecommand{\MR}{\relax\ifhmode\unskip\space\fi MR }

\providecommand{\href}[2]{#2}

\begin{document}
\section{Introduction}
Stein's method was introduced in the seminal paper \cite{stein1972bound}, and has since proved to be an extraordinarily useful tool in the area of normal approximation, see the text of \cite{chen2010normal} and references therein, and to an impressive, constantly expanding array of other distributions,  \cite{chen1975poisson}, \cite{luk1994stein}, \cite{pekoz1996stein}, \cite{chatterjee2011exponential},\cite{pekoz2013degree}, \cite{goldstein2018non}, \cite{gaunt2020stein}, as well as in concentration inequalities \cite{chatterjee2006stein}, \cite{chatterjee2010applications}, \cite{ghosh2011concentration},\cite{johnson2021concentration} stochastic analysis \cite{nourdin2012normal}, \cite{marinucci2015stein} \cite{last2016normal}, and data science \cite{erdogdu2016scaled}, \cite{liu2016kernelized}, \cite{liu2016stein}, \cite{fathi2020relaxing}. 

The method proceeds by the use of a characterizing equation for some target distribution which is to serve as an approximation to that of a random variable $W$ of interest. In the case of the normal, $W$ has the 
$\mathcal{N}(\mu,\sigma^2)$ distribution if and only if 
\begin{align}\label{eq:stein.char.general}
	\E[(W-\mu)f(W)]=\sigma^2 \E[f'(W)] \qmq{for all $f \in \mathcal{F}$,}
\end{align}
where here, and in like displays that follow in this section, we implicitly take $\mathcal{F}$ to be the class of functions for which the quantities written exist, which in particular will always include the collection of infinitely differentiable functions with compact support. 

Next, given any collection of functions $\mathcal{H}$, we define the pseudo-metric 
\begin{align}\label{pseudo-metric}
	d_{\mathcal H}(X,Y) = \sup_{h \in \mathcal{H}}|\E h(Y) - \E h(X)|,
\end{align}
which, for example, produces the Wasserstein $d(\cdot,\cdot)$ and Kolmogorov metric $d_K(\cdot,\cdot)$ by taking $\mathcal{H}$ to be the class 
$$
{\rm Lip}_1 = \{h: |h(y)-h(x)| \le |y-x|, \{x,y\} \subset \mathbb{R}\}\,\,
\mbox{or}\,\,
\textrm {Ind}= \{h: h(x)={\bm 1}_{(-\infty,z]}(x), z \in \mathbb{R}\},
$$
of Lipschitz-1 functions, or
the collection of all indicators of right closed,  left infinite intervals, respectively. 

We briefly summarize how Stein's method achieves a normal approximation of a random variable $W$ with mean $\mu$ and finite, non zero variance $\sigma^2$. Motivated by the characterizing equation \eqref{eq:stein.char.general}, given $h$ from the class of test functions $\mathcal{H}$,  with $Z \sim \mathcal{N}(\mu,\sigma^2)$ and $Nh=\E[h(Z)]$,  let $f_h$ be the unique bounded solution to the Stein equation
\begin{multline} \label{eq:stein_eq}
	\sigma^2 f'(w)-(w-\mu)f(w)=h(w)-Nh \qmq{thus yielding}\\ \E[\sigma^2 f_h'(W)-(W-\mu)f_h(W)]=\E[h(W)-Nh]. 
\end{multline}
Hence, one can obtain bounds for the quantity on the right hand side of \eqref{pseudo-metric} in terms of the left hand side of \eqref{eq:stein_eq}, to wit,
\begin{align}\label{pim.in.F}
	d_{\mathcal H}(W,Z) = \sup_{h \in \mathcal{H}}|\E[\sigma^2 f_h'(W)-(W-\mu)f_h(W)]|.
\end{align}
To put this expression to real use requires bounds on the solution $f_h$, such as on the magnitude of its derivatives. Though the form \eqref{pim.in.F} may seem at first glance to be more complicated and more difficult to handle than \eqref{pseudo-metric}, the latter involves the distribution of only a single random variable. 

It is perhaps for this reason that one of the best known ways to handle the right hand side of \eqref{pim.in.F} is through coupling constructions. An advance in understanding coupling methods was achieved in \cite{chen2010stein} by the introduction of the family of Stein couplings, of which many previously known couplings, such as the exchangeable pair, and size bias couplings, are special cases. 

We say the triple of random variables $(W'',W,G)$ is a Stein coupling when 
\begin{align} \label{def:stein.coupling}
	\E[(W-\mu)f(W)] = \E[G(f(W'')-f(W))].
\end{align}


\begin{example}\label{ex:ex_pair} When $\mu=0$, and $(W',W)$ are exchangeable and satisfy the linear regression condition
\begin{align} \label{eq:regression.condition}
	\E[W''|W]=(1-\lambda)W
\end{align}
for some $\lambda \in (0,1)$, then we obtain Stein's classical exchangeable pair that satisfies
\begin{align} \label{exp.pair.identity}
	\E[Wf(W)] = \frac{1}{2\lambda}\E[(W''-W)(f(W'')-f(W))],
\end{align}
that is, \eqref{def:stein.coupling} holds with $(W'',W,(W''-W)/2 \lambda)$, and hence is a Stein coupling. 
\end{example}

\begin{example}\label{eq:size_bias}
For an example of a different character, recalling that for a random variable $W \ge 0$ with positive mean $\mu$, we say $W^s$ has the $W$-size bias distribution when 
\begin{align} \label{def:size.bias}
	\E[Wf(W)] = \mu \E[f(W^s)],
\end{align}
a size bias coupling can be seen as a Stein coupling by 
taking the triple to be $(W'',W,G)=(W^s,W,\mu)$ when the sized biased variable $W^s$ is defined on the same space as $W$. 
\end{example}
An inventive coupling distinct from the ones just described is Coupling 2A in  \cite{chen2010stein}, all together demonstrating the large range of possibilities covered. 

One key step in computing bounds in the Kolmogorov or Wasserstein metric
to the normal using Stein couplings typically requires bounding a variance of a conditional expectation; 
for the exchangeable pair and size bias couplings described above, one needs to respectively bound
\begin{align} \label{eq:that.dreaded.var.of.cond.exp}
	{\rm Var}\left\{ \E\left( (W''-W)^2|W\right)\right\} \qmq{and} {\rm Var}\left\{ \E\left( W''-W|W\right) \right\};
\end{align}
see, for instance, the terms $r_1,r_2$ and $r_3$ in Theorem 2.1 of \cite{chen2010stein} for the Kolmogorov bound, or Theorem 3.20 of \cite{ross2011fundamentals}, for the Wasserstein. However, these terms are often difficult, tedious or even not tractable to compute. Bounds using size bias in  \cite{goldstein2011berry} for the lightbulb problem, described in detail below,   required a delicate and highly detailed analysis involving the eigenvaules of a certain Markov chain of multinomial type.  

On the other hand, the zero bias coupling \cite{goldstein1997stein} (see also \cite{chen2010normal}) can yield Wasserstein and Kolmogorov bounds without requiring cumbersome computation. We recall that for a mean $\mu$ random variable $W$ with variance $\sigma^2 \in (0,\infty)$, we say $W^*$ has the $W$-zero bias distribution on a space with some measure $Q$,  when 
\begin{align} \label{def:zero.bias}
	\E[(W-\mu)f(W)]=\sigma^2 \E_Q[f'(W^*)],
\end{align}
where we have subtracted the mean as in \cite{dobler2015new}, and consequently have $W^*=_d(W-\mu)^*+\mu$. In the canonical case where $\mu=0$ and $\sigma^2=1$, we obtain the bound (see \cite{goldstein2004normal}, \cite{goldstein2007l1})
$$
d(Z,W) \le 2d(W,W^*) 
$$
in the Wasserstein distance $d$ between $W$ and the standard normal $Z$. The right hand side above can be upper bounded by $2\E|W^*-W|$ for any coupling of $(W^*,W)$, 
in contrast to \eqref{eq:that.dreaded.var.of.cond.exp}. Kolmogorov bounds to the normal can be obtained without computation from any almost sure bound on $|W^*-W|$, see \cite{chen2010normal}. 

Theorem \ref{thm:normal.bounds} below provides bounds that generalize results previously obtained in both these metrics in the framework of \ycolor{Z}ero \ycolor{b}iased \ycolor{e}nhanced \ycolor{St}ein couplings, or \ycolor{zbest} for short, that includes both Stein and zero bias couplings as special cases.
For a random variable $W$ with mean $\mu$ and variance $\sigma^2 \in (0,\infty)$ on a probability space with measure $P$ (implicit in the expectation $\E$), we say the triple
$(W'',W',G)$ of random variables is a zbest coupling with respect to a measure $Q$ when 
\begin{align}\label{eq:it.is.zbest}
	\E[(W-\mu)f(W)] = \E_{Q}[G(f(W'')-f(W'))].
\end{align}
Clearly, zbest couplings include Stein couplings as a special case, though in contrast to \eqref{def:stein.coupling} $W'$ is not constrained to be equal to, nor even have the same distribution as, the random variable $W$. Enlarging the set of possibilities in this way leads to some additional useful flexibility, including generalizations of zero bias that  
likewise avoid the computation of quantities such as \eqref{eq:that.dreaded.var.of.cond.exp}.  

Naturally, nothing being for free, the work that would otherwise be required for the computation of \eqref{eq:that.dreaded.var.of.cond.exp} must be paid for in another manner.  In order to explain what is needed we give more detail on how a useful `second order' zbest coupling construction proceeds.  Lemma \ref{lem:zbest.lemma} in Section \ref{sec:zbest}  shows how we may begin with some initial zbest coupling, say, for instance, a standard exchangeable pair, and by changing the measure to some judicious $P^\dagger$,  be naturally led to a new zbest coupling that may possess favorable properties. However, this new coupling satisfies \eqref{eq:it.is.zbest} with $Q=P^\dagger$, and for Theorem \ref{thm:normal.bounds} to apply, it must be constructed on the original space.  That is, once the distribution of this new, advantageous 
triple is specified under $P^\dagger$,  we need to construct a triple $(W^\ddagger,W^\dagger,G^\dagger)$ on the original space so that ${\mathcal L}_{P}(W^\ddagger,W^\dagger,G^\dagger)= {\mathcal L}_{P^\dagger}(W'',W',G)$
where ${\mathcal L}$ denotes distribution, or law, on a space with underlying measure as indicated; the understanding of the procedure and the notation may be furthered upon reviewing Example \ref{ex:Bernolli}, where we treat an independent sum of Bernoulli variables. 

In general,  the final bound in Theorem \ref{thm:normal.bounds} depends on a number of choices. Indeed, this is so even when restricting to the use of Stein couplings, while for zbest one also has some influence when making the change of measure. Overall though, natural choices do seem to arise for given applications.

We apply the methods developed here for the lightbulb problem first analyzed in \cite{rao2007one} and further in \cite{goldstein2011berry}.
This instance demonstrates that the approach considered here may lead to easier computation of Kolmogorov and Wasserstein bounds which may also be superior to those produced by methods that require computation of quantities such as those in \eqref{eq:that.dreaded.var.of.cond.exp}.

In Section \ref{sec:zbest} we show how to use one zbest coupling 
to produce another by a change of 
measure,  obtain bounds in the Wasserstein and Kolmogorov metric, and give an introductory example. The Lightbulb Problem application can be found in Section \ref{sec:lightbulbs}.

\section{Zbest couplings} \label{sec:zbest}
Taking our variables as in \eqref{eq:it.is.zbest}, let
\begin{align} \label{def:D.lin_int}
	D=W''-W' \qmq{and}	W^*=UW'' + (1-U) W', 
\end{align}
where $U$ is a standard uniform variable independent of $W'', W'$ and $G$. For a smooth function $f$ for which the following expectations exist, and arbitrarily defining any $0/0$ expressions, integrating over the uniform variable for the last equality we obtain
\begin{align*} 
	G(f(W'')-f(W'))=
	DG\left[ \frac{f(W'')-f(W')}{W''-W'}\right]
	=  \E_Q[DG f'(W^*)|W'',W',G].
\end{align*}
Recalling \eqref{eq:it.is.zbest} and taking $Q$ expectation on both sides above results in the identity
\begin{align} \label{eq:z.b.and.GD}
	\E[(W-\mu)f(W)]\
	=\E_Q[G(f(W'')-f(W'))] = \E_Q[DGf'(W^*)].
\end{align}
Specializing to $f(w)=w$ in \eqref{eq:z.b.and.GD} we obtain
\begin{align} \label{eq:E[DG]=sigma2}
	\E_Q[DG]=\sigma^2, 
\end{align}
showing in particular that $Q(DG>0)>0$. 
When $DG=\sigma^2$ $Q$-a.s.\ then \eqref{eq:z.b.and.GD} recovers \eqref{def:zero.bias}, and hence occurs if and only if $W^*$ has the $W$-zero bias distribution.

The important message of Lemma \ref{lem:zbest.lemma} below is that when $D$ and $G$ in a given zbest coupling satisfy $DG \ge 0$ $Q$-a.s.\  then by  
biasing via the Radon Nikodym factor $DG/\E_Q[DG]$ to obtain the measure $P^\dagger$ one produces a new zbest coupling with $DG=\sigma^2$ $P^\dagger$-a.s, and therefore one that yields a variable with the $W$-zero bias distribution.  Further, if the non-negativity condition is not satisfied, one may nevertheless be able to construct a variable with a distribution close to that of the $W$-zero bias, and achieve a bound to the normal that has an additional term, see Theorem \ref{thm:normal.bounds}.

We illustrate the construction of the zero bias distribution as described in Lemma \ref{lem:zbest.lemma} in a few simple examples. This construction was previously understood only in the special case of the Stein exchangeable pair $(W'',W')$ as described in Example \ref{ex:ex_pair}, where it was shown that
\eqref{def:stein.coupling} holds with $\mu=0$ and $G=(W''-W')/2\lambda$. In particular, $DG=(W''-W')^2/2\lambda$ is non-negative, and non-trivial by \eqref{eq:E[DG]=sigma2}, hence one may construct a new measure $P^\dagger$ by
$$
dP^\dagger= \frac{(w''-w')^2}{2 \lambda \sigma^2}dP.
$$
Lemma \ref{lem:zbest.lemma} shows that this construction produces a new zbest coupling that yields a variable with the zero bias distribution via the interpolation \eqref{def:D.lin_int}.

A different, and rather simple special case is produced by taking a random variable $W$ with mean $\mu=0$, noting that $(W,0,W)$ trivially satisfies \eqref{eq:it.is.zbest} with $Q=P$ and that $DG=W^2$ is non-negative. In this case, as is known, and follows from Lemma \ref{lem:zbest.lemma}, by taking 
\begin{align} \label{eq:Wbox}
	dP^\square=\frac{w^2}{\sigma^2}dP \qmq{via \eqref{def:D.lin_int} one finds that} W^*=UW^\square
\end{align}
has the $W$-zero bias distribution, where $\mathcal{L}(W^\square)\sim \mathcal{L}_{P^\square}(W)$, and is independent of $U$.

\begin{lemma} \label{lem:zbest.lemma} Let \eqref{eq:it.is.zbest} hold, set
	$D=W''-W'$, let $R$ be a random variable satisfying $Q(R \ge 0)=1,\E_Q[R]=1$, and
	\begin{align}\label{eq:R=0.subset}
		\{D \not =0\}  \subseteq \{R\not =0\} \qm{a.s. $Q$.}
	\end{align}
	Then $(W'',W',G/R)$ is a zbest coupling with respect to the measure 
	\begin{align*} 
		dP^\dagger = R dQ .
	\end{align*}
	If
	\begin{align} \label{DG.Ind.ineq}
		DG \ge 0 \qmq{and} \{D\not =0\} \subseteq \{G\not =0\} \qmq{a.s.\ $Q$,}
	\end{align}
	then one may choose $R=DG/\sigma^2$, and the resulting $W^*$ of \eqref{def:D.lin_int} has the $W$-zero bias distribution under $P^\dagger$. 
\end{lemma}

\noindent {\em Proof:} 
As $P^\dagger(R=0)=0$ the ratio $G/R$ is well defined under $P^\dagger$, and changing measure we have 
\begin{multline*}
	\E_{P^\dagger}\left[\left( \frac{G}{R}\right)(f(W'')-f(W'))\right]
	= \E_{P^\dagger}\left[\left( \frac{G}{R}\right){\bf 1}(R \not =0)(f(W'')-f(W''))\right]
	\\
	=\E_Q\left[\left( \frac{G}{R}\right)R{\bf 1}(R \not = 0)(f(W'')-f(W'))\right] =\E_Q\left[ G{\bf 1}(R \not = 0)(f(W'')-f(W'))\right]
	\\ 
	=\E_Q\left[ G(f(W'')-f(W'))\right]= \E[(W-\mu)f(W)],
\end{multline*}
where we have used 
\eqref{eq:it.is.zbest} for the final equality, and for the penultimate, via  \eqref{eq:R=0.subset} 
and noting that $\{D\not =0\} \supseteq \{f(W'')-f(W')\not=0\}$, we have
\begin{multline*}
	{\bf 1}(R\not = 0)(f(W'')-f(W')) = {\bf 1}(R\not = 0, D \not = 0)(f(W'')-f(W')) \\
	= {\bf 1}(D \not = 0)(f(W'')-f(W'))=f(W'')-f(W')  \qmq{a.s.$Q$.}
\end{multline*}
Hence $(W'',W',G/R)$ is a zbest coupling under $P^\dagger$. 

For the final claim, by \eqref{DG.Ind.ineq} the variable $R=DG/\sigma^2$ is non-negative, integrates to 1 under $Q$ by \eqref{eq:E[DG]=sigma2}, and satisfies \eqref{eq:R=0.subset}, as 
\begin{align*}
	{\bf 1}(D \not = 0) = {\bf 1}(D \not = 0){\bf 1}(G \not = 0)= {\bf 1}(DG \not = 0).
\end{align*}
The first claim gives that $(W'',W',\sigma^2/D)$ is a zbest coupling under $P^\dagger$, and the last claim now follows from \eqref{eq:z.b.and.GD} with $Q=P^\dagger$.
\bbox

\begin{remark} \label{zb.from.sb}
	One important instance not previously noted in the literature that Lemma \ref{lem:zbest.lemma} brings to light begins with the size bias coupling $(W'',W',G)=(W^s,W,\mu)$ described in Example \ref{eq:size_bias}.
	The inclusion in \eqref{DG.Ind.ineq} is trivially satisfied and the product $DG = \mu(W^s-W)$ is non-negative when the coupling is monotone, that is, when $W^s \ge W$. In that case, taking $R=\mu(W^s-W)/\sigma^2=(W^s-W)/\E[W^s-W]$, the linear interpolation $W^*=UW''+(1-U)W'$ by a standard uniform variable $U$ independent of $(W'',W')$ under the change of measure
	\begin{align}\label{eq:dF.daggers.size.bias}
		dP^\dagger= \frac{w''-w'}{\E[W''-W']}dP
	\end{align}
	yields  $W^*$ with the $W$-zero bias distribution as in \eqref{def:zero.bias}.
\end{remark}

For normal approximation, the following theorem yields bounds that imply existing results in \cite{goldstein2004normal} and \cite{chen2010normal} for zero-bias couplings in both the Wasserstein and Kolmogorov metrics, as in that case the term $\E|1-GD|$ vanishes; the additional generality is gained by not requiring that $W^*$ has the $W$-zero bias distribution exactly.

\begin{theorem} \label{thm:normal.bounds}
	Let $W$ be a random variable with mean zero and variance 1. Then for any coupling of $W$ and $W^*=UW'' + (1-U) W'$ from a zbest triple $(W'',W',G)$ independent of $U \sim \mathcal{U}[0,1]$, with $D=W''-W'$
	\begin{align} \label{eq:thm.wass}
		d(W,Z) \le 2 \E|W^*-W| + \sqrt{2/\pi}\E|1-E[G D|W^*]|,
	\end{align}
	and when $|W^*-W| \le \delta$ for some $\delta \ge 0$,
	\begin{align} \label{eq:thm.gen.ED}
		d_K(W,Z) 
		\le  \left(\frac{1}{\sqrt{2\pi}}
		+\frac{\sqrt{2\pi}}{4}+1\right)\delta+E|1-E[G D|W^*]|
	\end{align}
	where the constant multiplying $\delta$ in the bound is less than $2.03$.
\end{theorem}

\noindent {\em Proof}: For the Wasserstein bound take  $h \in {\rm Lip}_1$ and let $f_h$ be the unique bounded solution of the Stein equation \eqref{eq:stein_eq}. By \eqref{eq:stein_eq} and \eqref{eq:z.b.and.GD} with $Q=P$,
\begin{multline*}
	\E[h(W)-Nh]= \E[f_h'(W)-Wf_h(W)] = \E[f_h'(W)-G D f_h'(W^*)]\\
	= \E[f_h'(W)-f_h'(W^*)+ (1-G D) f_h'(W^*)] = \E[f_h'(W)-f_h'(W^*)+ f_h'(W^*)E[1-G D |W^*] ]\\
	\le 2\E|W^*-W| + \sqrt{2/\pi}\E|1-E[G D|W^*]|,
\end{multline*}
where we have applied the mean value theorem to obtain the first term, and the bounds $\|f_h''\| \le 2$ and $\|f_h'\| \le \sqrt{2/\pi}$ from (2.13) of Lemma 2.4 in \cite{chen2010normal}. 
Taking supremum on the left hand side over $h \in {\rm Lip}_1$ yields the result. 

For the Kolmogorov bound, given an arbitrary $z \in \mathbb{R}$, by Lemma 2.2 in \cite{chen2010normal},
the unique bounded solution $f_z$ of the Stein equation \eqref{eq:stein_eq} with $h(w)={\bm 1}_{(-\infty,z]}(w)$ is given by 
\begin{align}\label{eq:f_z}
	f_z(w) = \sqrt{2 \pi}e^{w^2/2} \left( 
	\Phi(w)(1-\Phi(z)){\bm 1}_{w \le z}+\Phi(z)(1-\Phi(w)){\bm 1}_{w>z}
	\right),
\end{align}
where $\Phi$ is the cumulative distribution function of the standard normal; we set $f_z'(z)=\lim_{w \uparrow z}f_z'(w)$ so that the first equality in \eqref{eq:stein_eq} holds for all $w \in \mathbb{R}$. Hence, 
\begin{multline*}
	f_{z-\delta}'(W^*)={\bf 1}(W^* \le z-\delta)-\Phi(z-\delta) + W^*f_{z-\delta}(W^*)\\
	\le {\bf 1}(W \le z)-\Phi(z-\delta) + W^*f_{z-\delta}(W^*),
\end{multline*}
and so
\begin{multline*}
	{\bf 1}(W \le z) - \Phi(z) = (\Phi(z-\delta)-\Phi(z))+{\bf 1}(W \le z) - \Phi(z-\delta)\\
	\ge -\frac{\delta}{\sqrt{2\pi}}+f_{z-\delta}'(W^*)-W^*f_{z-\delta}(W^*).
\end{multline*}
Taking expectation, applying \eqref{eq:z.b.and.GD} and using that $|f_z'(w)| \le 1$ for all real $w,z$ via (2.8) of Lemma 2.3 of \cite{chen2010normal} in the final inequality, and handling the $1-GD$ term as above, yields
\begin{multline} \label{eq:lower.bd}
	P(W \le z) - \Phi(z) 
	\ge -\frac{\delta}{\sqrt{2\pi}}+\E[f_{z-\delta}'(W^*)-W^*f_{z-\delta}(W^*)] \\
	= -\frac{\delta}{\sqrt{2\pi}}+\E[G D f_{z-\delta}'(W^*)-W^*f_{z-\delta}(W^*)+(1-G D)f_{z-\delta}'(W^*)] \\
	\ge -\frac{\delta}{\sqrt{2\pi}}-\E|Wf_{z-\delta}(W)-W^*f_{z-\delta}(W^*)|-\E|1-E[G D|W^*]|.
\end{multline}
For the second term, applying the bound
$$
|(w+u)f_z(w+u)-(w+v)f_z(w+v)| \le \left(|w|+\frac{\sqrt{2\pi}}{4}\right)(|u|+|v|) \qm{for all $u,v,w \in \mathbb{R}$}
$$
from (2.10) of Lemma 2.3 in \cite{chen2010normal}, and letting $\Delta=W^*-W$, we have
\begin{multline*}
	\E|f_{z-\delta}(W)-W^*f_{z-\delta}(W^*)| = \E|f_{z-\delta}(W)-(W+\Delta)f_{z-\delta}(W+\Delta)|\\
	\le  \E\left[\left(|W|+\frac{\sqrt{2\pi}}{4}\right)|\Delta|\right]\le \left(1+\frac{\sqrt{2\pi}}{4}\right)\delta,
\end{multline*}
using $\E |W| \le \sqrt{\E W^2}=1$.
Substituting into \eqref{eq:lower.bd} produces a lower bound on $P(W \le z) - \Phi(z)$. 
Proving an analogous upper bound by the same reasoning, and then taking supremum over $z \in \mathbb{R}$ yields the claim \eqref{eq:thm.gen.ED}. \bbox

\begin{example} \label{ex:Bernolli} We end this section with an illustration for the sum of independent Bernoulli (indicator) variables. For $i \in [n]:=\{1,\ldots,n\}$, let $X_i \in \{0,1\}$ take the value $1$ with probability $p_i \in (0,1)$. Let ${\bm X}=(X_1,\ldots,X_n)$ and $Y=\sum_{i=1}^n X_i$, which has mean $\mu=\sum_{i=1}^n p_i$ and (positive) variance $\sigma^2 = \sum_{i=1}^n p_i(1-p_i)$. With $W=(Y-\mu)/\sigma$ and $Z$ the standard normal, we here obtain the bounds
\begin{align} \label{eq:bern.bounds}
d(W,Z) \le \frac{1}{\sigma} \qmq{and} d_K(W,Z) \le \frac{2.03}{\sigma}
\end{align}
using the size bias coupling discussed in Remark \ref{zb.from.sb} and Theorem \ref{thm:normal.bounds}. This case serves as a model for our approach in Section \ref{sec:lightbulbs}. 

In view of \eqref{def:size.bias}, as for all continuous functions $f$ and $i \in [n]$ we have
$$\E[X_i]\E[f(X_i^s)]=\E[X_i f(X_i)]= p_i f(1) = \E[X_i]\E[f(1)],$$
and see that size biasing any non-trivial indicator yields a variable that is identically equal to one. Further, when $Y$ is the sum of indicators, by Lemma 2.1 of \cite{goldstein1996multivariate}, a variable $Y^s$ with the $Y$-sized biased distribution may be constructed by 
sampling an index $I \in [n]$, independent of $\{X_i,i \in [n]\}$, with $P(I=i)=\E[X_i]/\mu$, setting $Y^s=Y$ if $X_I=1$ and otherwise, when $I=i$, letting $Y^s=\sum_{j=1}^n X_j^i$ where ${\bm X}^i=(X_1^i,\ldots,X_n^i)$ satisfies
\begin{align} \label{cond.dist.inds}
	P({\bm X}^i={\bm e})=P({\bm X}={\bm e}|X_i=1) \qmq{for all ${\bm e} \in \{0,1\}^n$.}
\end{align}
That is, in this second case one sets the $I^{th}$ indicator $X_I$ equal to 1 and constructs the remaining variables according to their conditional distribution given that updated value. 

As the Bernoullis here are independent, and independent of $I$, for $i\in [n]$ we have 
\begin{multline} \label{eq:samp.X.I}
P({\bm X}={\bm e},I=i)=P({\bm X}={\bm e})P(I=i) \\ \mbox{where} \quad P({\bm X}={\bm e}) = \prod_{i=1}^n p_i^{e_i}(1-p_i)^{1-e_i} \qmq{and} P(I=i)= \frac{p_i}{\mu}.
\end{multline}
We can construct a vector of indicators ${\bm X}''$ on the same space as ${\bm X}$ and $I$ according to the recipe suggested by \eqref{cond.dist.inds},  whose coordinate sum will be the $Y$-size biased variable $Y''$,  by specifying that
\begin{align} \label{eq:Bernoulli.X''_X'}
P({\bm X}''={\bm e}'',{\bm X}={\bm e},I=i)= P({\bm X}={\bm e})P(I=i){\bm 1}(e_j''=e_j, j \not =i, e_i''=1),
\end{align}
as, again  by independence, the distribution of $\{X_j, j \not =i\}$ is unaffected by conditioning on $X_i$.  Hence, the coupling amounts to simply replacing $X_I$ by 1 no matter its previous value when forming ${\bf X}''$, and when doing so we obtain 
\begin{align} \label{eq:1-X_I}
Y''-Y= (Y-X_I+1)-Y=1-X_I.
\end{align}
As $Y'' \ge Y$ the coupling is monotone, and following Remark \ref{zb.from.sb} we form the $P^\dagger$ distribution in \eqref{eq:dF.daggers.size.bias} by changing the $P$ measure in \eqref{eq:Bernoulli.X''_X'} according to \eqref{eq:1-X_I}, yielding
\begin{multline*}
	P^\dagger({\bm X}''={\bm e}'',{\bm X}'={\bm e}',I=i)=\frac{1-e_i}{\E[1-X_i]}P({\bm X}''={\bm e}'',{\bm X}'={\bm e}',I=i)\\
	=P(I=i){\bm 1}(e_i'=0,e_i''=1)\frac{P({\bm X}'={\bm e}')}{1-p_i}{\bm 1}(e_j''=e_j', j \not =i)\\
	=P(I=i){\bm 1}(e_i'=0,e_i''=1)\prod_{j \not =i} [P(X_j'=e_j'){\bm 1}(e_j''=e_j')].
\end{multline*}
Hence, sampling $({\bm X}, I)$ by \eqref{eq:samp.X.I} and setting $X_j^\ddagger=X_j^\dagger=X_j, j \not = I$, $X_I^\dagger=0,X_I^\ddagger=1$ yields
a coupling of $({\bm X}^\ddagger,{\bm X}^\dagger)$ and  ${\bm X}$ such that ${\mathcal L}_{P}({\bm X}^\ddagger,{\bm X}^\dagger)= {\mathcal L}_{P^\dagger}({\bm X}'',{\bm X}')$, under which
$$
Y^\dagger=Y-X_I \qmq{and} Y^\ddagger=Y-X_I+1 .
$$
Now by Lemma \ref{lem:zbest.lemma}, with $Y^*$ having the $Y$-zero bias distribution, we obtain
$$
Y^*-Y=UY^\ddagger+(1-U)Y^\dagger-Y=U(Y-X_I+1)+(1-U)(Y-X_I) -Y = U-X_I.
$$
Recalling that $U$ is independent of all other variables, taking expectation yields
$$
\E|Y^*-Y| = \E[(1-U)P(X_I=1)+UP(X_I=0)]=\frac{1}{2}, \qmq{and also} |Y^*-Y| \le 1. 
$$
The inequalities in \eqref{eq:bern.bounds} now follow by noting that the final terms of the bounds of Theorem \ref{thm:normal.bounds} vanish in the zero-bias case, and that $(aY)^*=aY^*$ for all $a \not =0$ (see \cite{chen2010normal}).

\end{example}

\section{The Lightbulb Process}\label{sec:lightbulbs}
The lightbulb problem was first considered in \cite{rao2007one} as a model for the behavior of skin receptors subject to a medication released by dermal patches, and subsequently studied in \cite{goldstein2011berry}. In the model, in a sequence of $n$ stages, $n$ skin receptors, which here we will imagine as lightbulbs, are toggled from one state to the other upon absorbing a pharmaceutical. Initially all $n$ lightbulbs are turned off. Then, at stages $r \in [n]$, a set of $r$ lightbulbs, selected uniformly from all subsets of $[n]$ of size $r$, independently of previous stages, have their state toggled.  
The random variable $Y$ of interest in the pharmaceutical application counts the number of lightbulbs that are turned on after stage $n$ is complete. The restriction that the number of stages and bulbs are equal is not essential, and more general frameworks are discussed in  \cite{rao2007one} and \cite{goldstein2011berry}, such as where in stage $r$ the number of toggled lightbulbs is some given $s_r \in [n]$; these variations can be analyzed as below.

More formally consider ${\bm X} \in \{0,1\}^{n \times n}$, a matrix of 
Bernoulli variables, called here a configuration, and whose components we refer to as toggle variables. The initial state of the bulbs is given deterministically with all bulbs off.
For stages $r \in [n]$ the components of ${\bm X}$
have the interpretation that 
\beas 
X_{ri} &=& \left\{
\begin{array}{cc}
	1 & \mbox{ if the status of bulb $i$ is changed at stage $r$,} \\
	0 &\mbox{ otherwise.}
\end{array}
\right.
\enas
As exactly $r$ of the $n$ bulbs are chosen uniformly to be toggled at stage $r$, taking
\begin{align*} 
	\mathcal{E}= \left\{ {\bm e} \in \{0,1\}^{n \times n}: \sum_{i=1}^n e_{ri}=r, r=1,\ldots,n \right\}, 
\end{align*}
and the toggles in each stage $r \in [n]$  independent, the distribution of ${\bm X}$ is given by
\bea \label{Xik-distribution}
P({\bm X}={\bf e})
=\left\{
\begin{array}{cl}
	\prod_{r=1}^n {n \choose r}^{-1} & \mbox{${\bm e} \in \mathcal{E}$}  \\
	0 & \mbox{otherwise.}
\end{array}
\right.
\ena
The toggle variables at stage $r$ that form the vector $(X_{r1},\ldots,X_{rn})$
are clearly exchangeable and the marginal distribution of the components $X_{ri}$ are Bernoulli with success probability $r/n$. For bulbs $i \in [n]$,
the variables
\bea
\label{def-Yk}
X_i = \left(\sum_{r=1}^n X_{ri}\right) \mbox{ mod }2 \qmq{and} Y=\sum_{i=1}^n X_i
\ena
are the indicator that
bulb $i$ is on at the terminal time $n$, and the total number of bulbs on at that time, respectively. In what follows we let, say, $Y'$ and $Y''$ be computed from configurations ${\bm X}'$ and ${\bm X}''$ as $Y$ is from ${\bm X}$.

In  \cite{goldstein2011berry} a normal approximation to $Y$ was achieved by constructing $Y''$ on the same space as $Y$ with the $Y$-size bias distribution, and using that construction \cite{ghosh2011concentration} also obtains a concentration bound. After reviewing that construction we apply Lemma \ref{lem:zbest.lemma} to create another zbest coupling in a way similar to the process followed in Example \ref{ex:Bernolli}.

Given a configuration ${\bm X} 
\in \mathcal{E}$, for a stage $r \in [n]$ and indices $\{i,j\} \subset [n]$, let ${\bm X}^{r, i \leftrightarrow j}$ be the configuration with components 
\beas
X_{sk}^{r, i \leftrightarrow j}=\left\{
\begin{array}{cl}
	X_{sk} & s \not =r\\
	X_{sk}& s=r, k \not \in \{i,j\}\\
	X_{sj} & s=r, k=i\\
	X_{si}& s=r, k=j,
\end{array}
\right.
\enas
that is, the new configuration is the same as ${\bm X}$, but with the toggle variables $X_{ri}$ and $X_{rj}$ interchanged. 

We take $n \ge 4$ and even for simplicity; Remark \ref{rem:odd} describes how the odd case may be handled. The size biased coupling $(Y'',Y)$ was constructed in \cite{goldstein2011berry} as follows. First, sample ${\bm X}$ according to \eqref{Xik-distribution}. Given ${\bm X}$, sample $I$ uniformly from $[n]$, and given ${\bm X}$ and $I$, let $J$ be an index chosen uniformly at random over the set of $n/2$ indices $\{j \in [n]:X_{n/2,j} \not = X_{n/2,I}\}$. Finally, let  ${\bm X}''={\bm X}$ if $X_I=1$ and ${\bm X}^{n/2, I \leftrightarrow J}$ otherwise.
Formally, with $e_i = (\sum_r e_{ri})\, {\rm mod}\,2$  and ${\bm e}, {\bm e}''\in \mathcal{E}$, we take
\begin{multline} \label{eq:dist.XIJ}
	P({\bm X}''={\bm e}'',{\bm X}= {\bm e}, I=i, J=j) \\
	= \frac{2{\bf 1}(e_{n/2,i} \not = e_{n/2,j}) }{n^2 \prod_{k=1}^n {n \choose r}}\left( {\bf 1}(e_i=1,{\bm e}''={\bm e})+{\bf 1}(e_i=0,{\bm e}''={\bm e}^{n/2, i \leftrightarrow j}) \right).
\end{multline}
Taking marginals, one easily finds that
\begin{align} \label{def:eta}
	P({\bm X}={\bm e}, I=i, J=j) = \eta {\bf 1}(e_{n/2,i} \not = e_{n/2,j}) \qmq{where}
	\eta= \frac{2}{n^2 \prod_{k=1}^n {n \choose r}}.
\end{align}

In \cite{goldstein2011berry}, via a verification of \eqref{cond.dist.inds}, it was shown that $Y''$, the number of lightbulbs that are on in the terminal stage in configuration ${\bm X}''$, has the $Y$ size bias distribution and satisfies
\begin{align} \label{eq:Ws.min.W}
	Y''- Y = 2{\bm 1}\{X_I=0,X_J=0\} = 2{\bm 1}\{Y'' \not = Y\}.
\end{align}
Referring the reader to \cite{goldstein2011berry} for a proof of the size property, we show \eqref{eq:Ws.min.W}. 
Indeed, if $X_I=1$ then $Y''=Y$ and the three quantities above are all zero. If $X_I=0$ and $X_J=1$ then all quantities are again zero as interchanging the toggles of bulbs $I$ and $J$ in stage $n/2$ will result in a configuration ${\bm X}''$ in which $X_I''=1, X_J''=0$, thus leaving $Y''$ with the same value as $Y$. In the remaining case $X_I=0,X_J=0$ and all quantities above equal 2, as the exchange toggles the final state of both bulbs $I$ and $J$, turning both from off to on, increasing $Y$ by 2.

Identity \eqref{eq:Ws.min.W} demonstrates that this size bias coupling is monotone, and thus Remark \ref{zb.from.sb} is in force. To continue, one needs to construct the biased version of the original pair according to $P^\dagger$ in \eqref{eq:dF.daggers.size.bias}, which, making use 
of the first equality in \eqref{eq:Ws.min.W},  is seen to be given by 
\begin{multline} \label{eq:joint.disc.daggers}
	P^\dagger({\bf X}''={\bm e}'',{\bm X}'={\bm e}',I=i,J=j)=\frac{y''-y'}{\E[Y''-Y']}P({\bm X}''={\bm e}'',{\bm X}'= {\bm e}', I=i, J=j)\\=
	\frac{{\bf 1}(e_i'=0,e_j'=0)}{P(Y'' \not = Y')}P({\bm X}''={\bm e}'',{\bm X}'= {\bm e}', I=i, J=j)\\
	= \frac{\eta}{P(Y'' \not = Y')}{\bf 1}(e_{n/2,i}' \not = e_{n/2,j}', e_i'=0,e_j'=0, {\bm e}''=({\bm e}')^{n/2, i \leftrightarrow j}),
\end{multline}
where in the last equality we have applied \eqref{eq:dist.XIJ} and the definition of $\eta$ in \eqref{def:eta}.
Integrating out the distribution of ${\bf X}''$, we find the marginal of ${\bm X}'$ under $P^\dagger$ satisfies
\begin{align} \label{eq:support.Xdagger.uniform}
	P^\dagger({\bm X}'={\bm e}', I=i, J=j)
	=\frac{\eta}{P(Y'' \not = Y')} {\bf 1}(e_{n/2,i}' \not = e_{n/2,j}', e_i'=0,e_j'=0).
\end{align}

\begin{lemma} \label{lem:exists.bulb.zcoupling}
	There exists a coupling of ${\bm X}$ and $({\bm X}^\ddagger,{\bm X^\dagger})$, having distributions specified by \eqref{def:eta} and \eqref{eq:joint.disc.daggers}, respectively, that satisfies
	\begin{align} \label{eq:coupling.diff.atmost2}
		|Y^\dagger-Y| \le 2 \qmq{and} Y^\ddagger = Y^\dagger + 2.
	\end{align}
\end{lemma}

\noindent {\em Proof:} To achieve the required construction, first sample ${\bm X},I,J$ from $P$ as in \eqref{def:eta}. Next, let $S$ be any random variable with support in $[n-1] \setminus \{n/2\}$, independent of $\{{\bm X},I,J\}$. Sample indices $K$ and $L$ that have marginal conditional laws $\mathcal{L}(\cdot|\cdot)$ specified by
\begin{multline*}
	\mathcal{L}(K|{\bm X},I,J,S) \sim \mathcal{U}\{k \not = J: X_{S,k} \not = X_{S,I}\} \qmq{and} \\
	\mathcal{L}(L|{\bm X},I,J,S) \sim \mathcal{U}\{l \not =I: X_{S,l} \not = X_{S,J}\}.
\end{multline*}
Now for $(a,b) \in \{0,1\}^2$, suppressing the dependence of $\phi_{ab}$ on $I,J,K,L$ and $S$, let
\begin{align*}
	{\bm X}^\dagger = \phi_{X_I,X_J}({\bm X}) \qmq{where} \phi_{ab}({\bm e})=\left\{
	\begin{array}{lc}
		\phi_{00}({\bm e}):={\bm e} & \mbox{if} \quad (a,b) = (0,0)\\
		\phi_{01}({\bm e}):={\bm e}^{S,J \leftrightarrow L} & \mbox{if} \quad (a,b) = (0,1)\\
		\phi_{10}({\bm e}):={\bm e}^{S,I \leftrightarrow K} & \mbox{if} \quad (a,b) = (1,0) \\
		\phi_{11}({\bm e}):={\bm e}^{n/2,I \leftrightarrow J} & \mbox{if} \quad (a,b)=(1,1).
	\end{array}
	\right.
\end{align*}

One easily verifies that $P({\bm X}^\dagger \in \mathcal{E})=1$ and that each $\phi_{ab}$ is an involution mapping between $\{{\bm e} \in \mathcal{E}: e_I=a,e_J=b\}$ and $\{{\bm e} \in \mathcal{E}:e_I=0,e_J=0\}$. 
Moreover, for any ${\bm f}\in \mathcal{E}$,
\begin{multline*}
	P({\bm X}^\dagger = {\bm f}, I=i,J=j) =  P(\phi_{X_i,X_j}({\bm X}) = {\bm f}, I=i,J=j)\\
	= \sum_{(a,b) \in \{0,1\}^2} P(\phi_{ab}({\bm X})= {\bm f},I=i,J=j, X_i=a,X_j=b).
	\\
	= \sum_{(a,b) \in \{0,1\}^2} P({\bm X}= \phi_{ab}({\bm f}),I=i,J=j, X_i=a,X_j=b).
\end{multline*}

Note that on the event $\{I=i,J=j\}$ we have $X_{n/2,i} \not = X_{n/2,j}$, hence the probability in the $a,b^{th}$ summand is zero unless
$(\phi_{ab}({\bm f}))_{n/2,i} \not = (\phi_{ab}({\bm f}))_{n/2,j}$, which is equivalent to the condition that $f_{n/2,i} \not = f_{n/2,j}$.  Similarly, for this probability to be non-zero we must have that $(\phi_{ab}({\bm f}))_i=a, (\phi_{ab}({\bm f}))_j=b$, which implies that $f_i=0, f_j=0$. 
For ${\bm f}\in\mathcal{E}$ such that $f_{n/2,i} \not = f_{n/2,j}$ and  $f_i=0,f_j=0$ the restriction in the summand that $X_i=a,X_j=b$ is redundant, and as $\phi_{ab}({\bm f}) \in \mathcal{E}$ and the probability $P({\bm X}={\bm e},I=i,J=j)$ is constant over its support by \eqref{def:eta}, we obtain
\begin{align*}
	P({\bm X}^\dagger = {\bm f}, I=i,J=j) \propto {\bm 1}(f_{n/2,i} \not = f_{n/2,j}, f_i=0,f_j=0),
\end{align*}
in  agreement with \eqref{eq:support.Xdagger.uniform}. 
That is, 
the constructed ${\bm X}^\dagger$ has the desired distribution, and setting ${\bm X}^\ddagger= \phi_{11}({\bm X}^\dagger)$ it is immediate that the pair $({\bm X}^\ddagger,{\bm X}^\dagger)$ has joint distribution \eqref{eq:joint.disc.daggers}.

Lastly, as each mapping $\phi_{ab}$ transposes at most two toggles, the final states of ${\bm X}$ and ${\bm X}^\dagger$  can differ for at most two lightbulbs, 
thus verifying the first claim of \eqref{eq:coupling.diff.atmost2}. The second claim follows directly from \eqref{eq:joint.disc.daggers}.
\bbox

The mean $\mu$ and variance $\sigma^2$ of $Y$ are given by
\begin{align*}
	\mu=n/2 \qmq{and}	\sigma^2= \frac{n}{4}\left( 
	1+(n-1)\lambda_n
	\right) \qmq{with} \lambda_n =\prod_{s=1}^n \left(1-\frac{4s(n-s)}{n(n-1)}  \right),
\end{align*}
where the first equality follows from a simple symmetry argument, and the next from Section 2.3 of \cite{rao2007one}, and as given in (3), (4) and (5) of \cite{goldstein2011berry}, where $\lambda_n$ there is denoted $\lambda_{n,2,{\bm n}}$. Note that when $n$ is even, as the terms in the product for $s$ and $n-s$ are equal, 
\begin{align*}
	\lambda_n = \left(1-\frac{4(n/2)^2}{n(n-1)}\right) \prod_{s=1}^{n/2-1}
	\left(1-\frac{4s(n-s)}{n(n-1)}  \right)^2 = -\frac{1}{n-1}\prod_{s=1}^{n/2-1}
	\left(1-\frac{4s(n-s)}{n(n-1)}  \right)^2 \le 0, 
\end{align*}
and hence $\sigma^2 \le n/4$. Applying the reasoning that proves (1.6) of 
\cite{englund1981remainder} demonstrates that the order $1/\sigma$ in the Kolmogorov bound in \cite{goldstein2011berry}, and \eqref{eq:8.12}, is unimprovable. 

\begin{theorem}
	Let $Y$ be the number of lightbulbs switched on in the final stage of the lightbulb process with $n \ge 4$ even stages. Then, with $\mu=\E[Y],\sigma^2={\rm Var}(Y)$ and $W=(Y-\mu)/\sigma$, 
	with $Z \sim \mathcal{N}(0,1)$, 
	\begin{align} \label{eq:8.12}
		d(W,Z) \le 6/\sigma \qmq{and}	d_K(W,Z) \le 8.12/\sigma.
	\end{align}
\end{theorem}
\noindent {\em Proof:} Letting $U \sim \mathcal{U}[0,1]$ be independent of $Y,Y^\dagger,Y^\ddagger$ of Lemma \ref{lem:exists.bulb.zcoupling}, by Lemma \ref{lem:zbest.lemma}, 
\begin{align*}
	Y^*=UY^\ddagger + (1-U)Y^\dagger
\end{align*}
has the $Y$-zero bias distribution. Using \eqref{eq:coupling.diff.atmost2} of Lemma \ref{lem:exists.bulb.zcoupling}, 
\begin{align*}
	|Y^*-Y| = |U(Y^\dagger+2) +(1-U)Y^\dagger -Y| = |Y^\dagger -Y +2U| \le (2 + 2U)\le 4,
\end{align*}
and $E|Y^*-Y| \le 3$, and the claims in \eqref{eq:8.12} follow from these inequalities and \eqref{eq:thm.wass} and \eqref{eq:thm.gen.ED}, respectively, and the scaling $(aY)^*=_daY^*$ for $a \not =0$ as in Example \ref{ex:Bernolli}.
\bbox

With
$$
\overline{\Delta}_0 = \frac{1}{2\sqrt{n}} + \frac{1}{2n} + \frac{1}{3}e^{-n/2}
$$
and $B_n$ denoting the constant in \cite{goldstein2011berry} that for $n \ge 6$ replaces 8.12 in \eqref{eq:8.12}, over that range, using $\sigma^2 \le n/4=:\overline{\sigma^2}$, we have
\begin{multline*}
	B_n = 
	\frac{n}{2\sigma}\overline{\Delta}_0 + 1.64\frac{n}{\sigma^2}+2
	\ge  
	\frac{n}{2\overline{\sigma}} \frac{1}{2\sqrt{n}} + 1.64\frac{n}{\overline{\sigma}^2}+2 = \frac{1}{2}+4(1.64)+2=9.06.
\end{multline*}
Hence the bound produced here is superior for all values over which the former was valid.

\begin{remark}\label{rem:odd}
	The odd case can be handled using randomization as in \cite{goldstein2011berry}. In particular, one constructs a configuration of toggle variables ${\bm V}$ giving rise to an intermediate variable $V$, close to $Y$ and having favorable symmetry properties for size biasing, by randomly adding, or removing, a toggle variable in the two middle stages $m$ and $m+1$, respectively. To size bias $V$, with $I$ the random index as in the even case, when $V_I=0$ a randomization is applied to the toggle variable of bulb $I$ in one of the two middle stages. If that does not succeed in changing the state (the event $F=0$),  then an interchange of the toggle variables $V_I$ and $V_J$, such as achieved in \eqref{eq:dist.XIJ}, is performed in the chosen middle stage. This gives rise to an additional term ${\bf 1}(V_I=0,F=1)$ for the difference  \eqref{eq:Ws.min.W} (see (46) of \cite{goldstein2011berry}), that here produces our change of measure. That term, and one accounting for the difference between $V$ and $Y$, produce only a small additional term in the final bound; note the difference between the even and odd case results in Theorem 1 of \cite{goldstein2011berry}.
\end{remark}



\begin{thebibliography}{10}
	
	\bibitem{chatterjee2006stein}
	Sourav Chatterjee, \emph{{S}tein's method for concentration inequalities},
	arXiv preprint math/0604352 (2006).
	
	\bibitem{chatterjee2010applications}
	Sourav Chatterjee and Partha~S Dey, \emph{Applications of {S}tein’s method
		for concentration inequalities}, The Annals of Probability \textbf{38}
	(2010), no.~6, 2443--2485.
	
	\bibitem{chatterjee2011exponential}
	Sourav Chatterjee, Jason Fulman, and Adrian R{\"o}llin, \emph{Exponential
		approximation by {S}tein’s method and spectral graph theory}, ALEA Lat. Am.
	J. Probab. Math. Stat \textbf{8} (2011), no.~1, 197--223.
	
	\bibitem{chen1975poisson}
	Louis~HY Chen, \emph{Poisson approximation for dependent trials}, The Annals of
	Probability \textbf{3} (1975), no.~3, 534--545.
	
	\bibitem{chen2010normal}
	Louis~HY Chen, Larry Goldstein, and Qi-Man Shao, \emph{Normal approximation by
		{S}tein’s method}, Springer Science \& Business Media, 2010.
	
	\bibitem{chen2010stein}
	Louis~HY Chen and Adrian R{\"o}llin, \emph{{S}tein couplings for normal
		approximation}, arXiv preprint arXiv:1003.6039 (2010).
	
	\bibitem{dobler2015new}
	Christian D{\"o}bler, \emph{New {B}erry-{E}sseen and {W}asserstein bounds in
		the {CLT} for non-randomly centered random sums by probabilistic methods},
	arXiv preprint arXiv:1504.05938 (2015).
	
	\bibitem{englund1981remainder}
	Gunnar Englund, \emph{A remainder term estimate for the normal approximation in
		classical occupancy}, The Annals of Probability (1981), 684--692.
	
	\bibitem{erdogdu2016scaled}
	Murat~A Erdogdu, Lee~H Dicker, and Mohsen Bayati, \emph{Scaled least squares
		estimator for {GLM}s in large-scale problems}, Advances in Neural Information
	Processing Systems \textbf{29} (2016).
	
	\bibitem{fathi2020relaxing}
	Max Fathi, Larry Goldstein, Gesine Reinert, and Adrien Saumard, \emph{Relaxing
		the {G}aussian assumption in shrinkage and {SURE} in high dimension}, arXiv
	preprint arXiv:2004.01378 (2020).
	
	\bibitem{gaunt2020stein}
	Robert~E Gaunt, \emph{{S}tein factors for variance-gamma approximation in the
		{W}asserstein and {K}olmogorov distances}, arXiv preprint arXiv:2008.06088
	(2020).
	
	\bibitem{ghosh2011concentration}
	Subhankar Ghosh and Larry Goldstein, \emph{Concentration of measures via
		size-biased couplings}, Probability {T}heory and {R}elated fields
	\textbf{149} (2011), no.~1, 271--278.
	
	\bibitem{goldstein2004normal}
	Larry Goldstein, \emph{Normal approximation for hierarchical structures}, The
	{A}nnals of {A}pplied {P}robability \textbf{14} (2004), no.~4, 1950--1969.
	
	\bibitem{goldstein2007l1}
	\bysame, \emph{${L}^1$ bounds in normal approximation}, The Annals of
	Probability \textbf{35} (2007), no.~5, 1888--1930.
	
	\bibitem{goldstein2018non}
	\bysame, \emph{Non-asymptotic distributional bounds for the {D}ickman
		approximation of the running time of the {Q}uickselect algorithm}, Electronic
	{J}ournal of {P}robability \textbf{23} (2018), 1--13.
	
	\bibitem{goldstein1997stein}
	Larry Goldstein and Gesine Reinert, \emph{{S}tein's method and the zero bias
		transformation with application to simple random sampling}, The {A}nnals of
	{A}pplied {P}robability \textbf{7} (1997), no.~4, 935--952.
	
	\bibitem{goldstein1996multivariate}
	Larry Goldstein and Yosef Rinott, \emph{Multivariate normal approximations by
		{S}tein's method and size bias couplings}, Journal of Applied Probability
	\textbf{33} (1996), no.~1, 1--17.
	
	\bibitem{goldstein2011berry}
	Larry Goldstein and Haimeng Zhang, \emph{A {B}erry-{E}sseen bound for the
		lightbulb process}, Advances in {A}pplied {P}robability \textbf{43} (2011),
	no.~3, 875--898.
	
	\bibitem{johnson2021concentration}
	Tobias Johnson and Erol Pek{\"o}z, \emph{Concentration inequalities from
		monotone couplings for graphs, walks, trees and branching processes}, arXiv
	preprint arXiv:2108.02101 (2021).
	
	\bibitem{last2016normal}
	G{\"u}nter Last, Giovanni Peccati, and Matthias Schulte, \emph{Normal
		approximation on {P}oisson spaces: {M}ehler’s formula, second order
		{P}oincar{\'e} inequalities and stabilization}, Probability {T}heory and
	{R}elated {F}ields \textbf{165} (2016), no.~3, 667--723.
	
	\bibitem{liu2016kernelized}
	Qiang Liu, Jason Lee, and Michael Jordan, \emph{A kernelized {S}tein
		discrepancy for goodness-of-fit tests}, International conference on machine
	learning, PMLR, 2016, pp.~276--284.
	
	\bibitem{liu2016stein}
	Qiang Liu and Dilin Wang, \emph{Stein variational gradient descent: A general
		purpose {B}ayesian inference algorithm}, Advances in {N}eural {I}nformation
	{P}rocessing {S}ystems \textbf{29} (2016).
	
	\bibitem{luk1994stein}
	Ho~Ming Luk, \emph{Stein's method for the gamma distribution and related
		statistical applications}, Ph.D. thesis, University of Southern California,
	1994.
	
	\bibitem{marinucci2015stein}
	Domenico Marinucci and Maurizia Rossi, \emph{Stein--{M}alliavin approximations
		for nonlinear functionals of random eigenfunctions on ${S}^d$}, Journal of
	{F}unctional {A}nalysis \textbf{268} (2015), no.~8, 2379--2420.
	
	\bibitem{nourdin2012normal}
	Ivan Nourdin and Giovanni Peccati, \emph{Normal approximations with {M}alliavin
		calculus: from {S}tein's method to universality}.
	
	\bibitem{pekoz1996stein}
	Erol~A Pek{\"o}z, \emph{Stein's method for geometric approximation}, Journal of
	{A}pplied {P}robability \textbf{33} (1996), no.~3, 707--713.
	
	\bibitem{pekoz2013degree}
	Erol~A Pek{\"o}z, Adrian R{\"o}llin, and Nathan Ross, \emph{Degree asymptotics
		with rates for preferential attachment random graphs}, The {A}nnals of
	{A}pplied {P}robability \textbf{23} (2013), no.~3, 1188--1218.
	
	\bibitem{rao2007one}
	C~Radhakrishna Rao, M~Bhaskara Rao, and Haimeng Zhang, \emph{One bulb? {T}wo
		bulbs? {H}ow many bulbs light up?—a discrete probability problem involving
		dermal patches}, Sankhy{\=a}: The {I}ndian {J}ournal of {S}tatistics (2007),
	137--161.
	
	\bibitem{ross2011fundamentals}
	Nathan Ross, \emph{Fundamentals of {S}tein’s method}, Probability {S}urveys
	\textbf{8} (2011), 210--293.
	
	\bibitem{stein1972bound}
	Charles Stein, \emph{A bound for the error in the normal approximation to the
		distribution of a sum of dependent random variables}, Proceedings of the
	sixth Berkeley symposium on mathematical statistics and probability, volume
	2: Probability theory, University of California Press, 1972, pp.~583--602.
	
\end{thebibliography}

%

\end{document}